\documentclass[amscd,amssymb,11pt,reqno]{amsart}

\newtheorem{Thm}{Theorem}[section]

\newtheorem{Cor}[Thm]{Corollary}
\newtheorem{Lem}[Thm]{Lemma}
\newtheorem{Prop}[Thm]{Proposition}
\newtheorem{Def}[Thm]{Definition}
\newtheorem{Rmk}[Thm]{Remark}

\usepackage{graphicx}
\usepackage{latexsym}
\begin{document}

%\begin{Large}

\vspace{1.5 cm}

\title[The lower dimensional hyperbolic Busemann-Petty  problem]
      {A solution to the lower dimensional\\ Busemann-Petty  problem  \\ in the hyperbolic space}

\author{  V.Yaskin}

\address{Department of Mathematics, University of Missouri, Columbia, MO 65211, USA.}

\email{yaskinv@math.missouri.edu}

%\begin{center}{\bf Sections} \end{center}

\begin{abstract}
The lower dimensional Busemann-Petty problem asks whe\-ther origin symmetric convex bodies in
$\mathbb{R}^n$ with smaller volume of all $k$-dimensional sections necessarily have smaller volume.
As proved by Bourgain and Zhang, the answer to this question is negative if $k>3$. The problem is
still open for $k=2,3$.  In this article we formulate and completely solve the lower dimensional
Busemann-Petty problem in the hyperbolic space $\mathbb{H}^n$.

\end{abstract}



\maketitle

\section{Introduction}

The Busemann-Petty problem  asks whether origin symmetric convex bodies in $\mathbb{R}^n$ with
smaller hyperplane sections necessarily have smaller volume. The answer to this problem is
affirmative if $n\le 4$ and negative if $n\ge 5$ (see \cite{GKS}, \cite{Zh2} or \cite[Chapter
5]{Kbook}  for the solution and historical details). In \cite{Y} the author solved the
Busemann-Petty problem in hyperbolic and spherical spaces.

The lower dimensional Busemann-Petty problem (LDBP) in $\mathbb{R}^n$ asks the same question with
$k$-dimensional subspaces in place of hyperplanes. Bourgain and Zhang \cite{BZ} proved that this
problem has a negative answer if $3<k<n$, see  \cite{K4} for another solution. The cases $k=2,3$
are still open in dimensions $n>4$.

In this paper we study the lower dimensional Busemann-Petty problem in the hyperbolic space.
Namely, let $1\le k<n$,  and  $K$, $L$ be origin-symmetric  convex bodies  in $\mathbb{H}^n$, $n\ge
3$, such that
$$\mathrm{vol}_{k}(K\cap H)\le \mathrm{vol}_{k}(L\cap H)$$ for every $k$-dimensional totally geodesic
plane through the origin. Does it follow that
$$\mathrm{vol}_{n}(K)\le  \mathrm{vol}_{n}(L)?$$

For the case $k=1$ the answer is trivially affirmative, since in all directions the radius of $K$
does not exceed the radius of $L$. In this paper we prove that the answer to the hyperbolic lower
dimensional Busemann-Petty problem is negative for every $2\le k <n$.

\section{Hyperbolic geometry}

It is well-known (see \cite[\S 10]{DFN} or \cite[\S 4.5]{R} )  that the hyperbolic space
$\mathbb{H}^n$ can be identified with the interior of the unit ball $B^n$ in $\mathbb{R}^{n}$ with
the metric:

\begin{eqnarray}\label{eqn:metric}
ds^2=4 \frac{dx_1^2+\cdots +dx_n^2}{(1- (x_1^2+\cdots+x_n^2))^{2}}.
\end{eqnarray}
This is called the Poincar\'{e} model of the hyperbolic space in the ball.   The geodesic lines in
this model are arcs of the circles orthogonal to the boundary of the ball $B^n$ and straight lines
through the origin.

Since for any two points in the hyperbolic space there exists a unique geodesic connecting them,
the definition of convexity in the hyperbolic  space will be analogous to that in the Euclidean
space (see \cite[Chapter I, \S 12]{P}).  A body $K$ (compact set with non-empty interior) is called
{\it convex}, if for every pair of points in $K$, the geodesic segment joining them also belongs to
the body $K$.

In order to distinguish between different types of convexity in the unit ball, we use the following
system of notations. Let $K$ be a body in the open unit ball $B^n$. The body $K$ is called
h-convex, if it is convex in the hyperbolic metric defined in the ball $B^n$. Similarly it is
called e-convex, if it is convex in the usual Euclidean sense. Analogously, h-geodesics are the
straight lines of the hyperbolic metric and e-geodesics are the usual Euclidean straight lines.

A submanifold $\mathcal F$ in a Riemannian space $\mathcal R$ is called {\it totally geodesic} if
every geodesic in $\mathcal F$ is also a geodesic in the space $\mathcal R$. In the Euclidean space
the totally geodesic submanifolds are Euclidean planes. %, on the sphere they are great subspheres.
In the  Poincar\'{e} model of the hyperbolic space described above the totally geodesic
submanifolds are represented by the spheres orthogonal to the boundary of the unit ball $B^n$ and
Euclidean planes through the origin. (We want to emphasize that a $k$-dimensional submanifold
passing through the origin is totally geodesic if and only if it is a $k$-dimensional Euclidean
plane).  In a sense, totally geodesic submanifolds are analogs of Euclidean planes in Riemannian
spaces. For elementary properties of totally geodesic submanifolds see \cite[Chap.5, \S 5]{A}.

The {\it Minkowski functional} of a star-shaped origin-symmetric body $K\subset \mathbb R^n$ is
defined as
$$\|x\|_K=\min \{a\ge 0: x \in aK \}.$$
The {\it radial function} of $K$ is given by $\rho_K(x) =\|x\|_K^{-1}$. If $x\in S^{n-1}$ then the
radial function $\rho_K(x)$ is the Euclidean distance from the origin to the boundary of $K$ in the
direction of $x$.

The volume element of the metric (\ref{eqn:metric}) equals

$$d\mu_n=2^n \frac{dx_1\cdots dx_n}{(1- (x_1^2+\cdots+x_n^2))^{n}}=2^{n} \frac{dx}{(1- |x|^2)^{n}}.$$
Therefore the hyperbolic volume of a body $K$ is given by the formula:
$$\mathrm{vol}_n(K)= \int_{K} d\mu_{n}=2^n \int_K \frac{dx}{(1- |x|^2)^{n}}.$$
Note that in the polar coordinates of $\mathbb{R}^n$ the latter formula looks as follows:
\begin{eqnarray}\label{eqn:polarvolume}
\mathrm{vol}_n(K)=2^n\int_{S^{n-1}}\int_{0}^{\|\theta\|^{-1}_K} \frac{r^{n-1}}{(1- r^2)^{n}}dr\
d\theta.
\end{eqnarray}
Similarly, if $H$ is a $k$-dimensional  hyperbolic totally geodesic plane through the origin (as
mentioned above, this is just a $k$-dimensional  Euclidean plane through the origin), then the
volume element of $H$ in the metric (\ref{eqn:metric}) is
$$d\mu_{k}=2^{k} \frac{dx}{(1- |x|^2)^{k}},$$
therefore the hyperbolic $k$-volume of the section of $K$ by $H$ is given by the formula:
$$\mathrm{vol}_k(K\cap H)=\int_{K\cap H} d\mu_{k}=2^{k} \int_{K\cap H} \frac{dx}{(1- |x|^2)^{k}},$$
or in polar coordinates:
\begin{eqnarray}\label{eqn:section-polarvolume}
\mathrm{vol}_k(K\cap H)=2^k\int_{S^{n-1}\cap H}\int_{0}^{\|\theta\|^{-1}_K} \frac{r^{k-1}}{(1-
r^2)^{k}}dr\ d\theta.
\end{eqnarray}

Even though our main object is hyperbolic geometry, let us briefly mention that,  along with the
hyperbolic and Euclidean metrics, we can define the spherical metric in the unit ball $B^n$:
$$ds^2=4 \frac{dx_1^2+\cdots +dx_n^2}{(1+(x_1^2+\cdots+x_n^2))^{2}}.$$
The geodesic lines in this model are arcs of the circles intersecting the boundary of the ball
$B^n$ in antipodal points and straight lines through the origin. Such lines will be called
s-geodesics. The body $K$ is called s-convex, if it is convex in the spherical metric defined in
the ball $B^n$. (This notion is well-defined, since in this model every two points can be joined by
a unique geodesic).

Finally, a simple observation about all introduced types of convexity is that any s-convex body
containing the origin is also e-convex and any e-convex body containing the origin is h-convex.
(See for example \cite{MP}).

\section{Fourier transform of distributions}

The Fourier transform of a distribution $f$ is defined by $\langle\hat{f}, {\phi}\rangle= \langle
f, \hat{\phi} \rangle$ for every test function $\phi$ from the Schwartz space $ \mathcal{S}$ of
rapidly decreasing infinitely differentiable functions on $\mathbb R^n$.

%A distribution is called {\it positive definite} if for every test function $\phi$
%$$\langle f, \phi \ast \overline{\phi(-x)}\rangle \ge 0.$$
%A distribution is positive definite if and only if its Fourier transform is a positive distribution
%(in the sense that $\langle \hat{f},\phi \rangle \ge 0$ for every non-negative test function
%$\phi$; see, for example, \cite[p.152]{GV}).

We say that a distribution $f$ is {\it positive definite} if its Fourier transform is a positive
distribution, in the sense that $\langle \hat{f},\phi \rangle \ge 0$ for every non-negative test
function $\phi$.

We say that a closed bounded set $K$ in $\mathbb{R}^n$ is a {\it star body} if for every $x\in K$
each point of the interval $[0,x)$ is an interior point of $K$, and  $\|x\|_K$, the Minkowski
functional of $K$, is a continuous function on $\mathbb{R}^n$.

Let $K$ be a star body and $\xi\in S^{n-1}$,  the {\it parallel section function} of $K$ is defined
as follows:
%$$A_{K,\xi}(z)=\int_{K\cap \{\langle x,\xi\rangle=z\}} dx.$$
$$A_{K,\xi}(z)=\mathrm{vol}_{n-1}({K\cap \{\langle x,\xi\rangle=z\}}).$$
(We also assume that $K\cap \{\langle x,\xi\rangle=z\}$ is star-shaped for small $z$). Recall the
following fact:

\begin{Thm}\label{Thm:GKS}{\rm (\cite{GKS}, Theorem 1)}
Let $K$ be an origin-symmetric star body in $\mathbb{R}^n$ with $C^\infty$ boundary, and let
$k\in\mathbb{N}\setminus \{0\}$, $k\ne n-1$. Suppose that $\xi\in S^{n-1}$, and let $A_\xi$ be the
corresponding parallel section function of $K$.
\newline
(a) If $k$ is even, then
$$(\|x\|_K^{-n+k+1})^\wedge(\xi)=(-1)^{k/2}\pi (n-k-1) A_\xi^{(k)}(0).$$
\newline
(b) If $k$ is odd, then
\begin{eqnarray*}
(\|x\|_K^{-n+k+1})^\wedge(\xi) =(-1)^{(k+1)/2}2(n-1-k)k!\times\hspace{3.5cm}\\
\times\int_0^\infty
\frac{A_\xi(z)-A_\xi(0)-{A''}_\xi(0)\frac{z^2}{2}-\cdots-A_\xi^{(k-1)}(0)\frac{z^{k-1}}{(k-1)!}}{z^{k+1}}dz,
\end{eqnarray*}
where $A_\xi^{(k)}$ stands for the derivative of the order $k$ and the Fourier transform is
considered in the sense of distributions.
\end{Thm}
In particular, it follows that for infinitely smooth bodies the Fourier transform of
$\|x\|^{-n+k+1}$ restricted to the unit sphere  is a continuous function (see also \cite[Section
3.2]{Kbook}). This remark explains why integration over the sphere in the next lemma makes sense.
The following is Parseval's formula on the sphere proved by  Koldobsky \cite{K3}.

\begin{Lem}\label{Lem:Parseval}
If $K$ and $L$ are origin symmetric infinitely smooth star bodies in $\mathbb{R}^n$ and $0<p<n$,
then

$$\int_{S^{n-1}} \left(\|x\|_K^{-p}\right)^\wedge (\xi) \left(\|x\|_L^{-n+p}\right)^\wedge
(\xi)d\xi= (2\pi)^n \int_{S^{n-1}} \|x\|_K^{-p}\|x\|_L^{-n+p} dx.$$
\end{Lem}

The following result  was also proved in \cite{K3}.

\begin{Lem}\label{Lem:FTonH}
Let $L$ be an origin symmetric star body with $C^{\infty}$ boundary in $\mathbb{R}^n$. Then for
every $(n-k)$-dimensional subspace $H$ of $\mathbb{R}^n$ we have
$$(2\pi)^k\int_{S^{n-1}\cap H} \|\theta\|_L^{-n+k} d\theta = \int_{S^{n-1}\cap H^\perp} (\|x\|_L^{-n+k})^\wedge(\theta)
d\theta.$$
\end{Lem}
The preceding two lemmas were formulated for Minkowski functionals, but in fact they are true for
arbitrary infinitely differentiable even functions on the sphere extended to
$\mathbb{R}^n\setminus\{0\}$ as homogeneous functions of corresponding degrees. (Indeed, any such
function of degree $-p$ can be obtained as the difference of Minkowski functionals raised to the
power $-p$).

The next lemma is a Fourier analytic version of a result of Zhang  \cite[Lemma 2]{Zh1}.

\begin{Lem}\label{Lem:Zhang} Let $k$ be an integer, $1\le k \le n-1$,
and let $f$ be an infinitely differentiable even function on the sphere ${S}^{n-1}$, such that
$f(x/|x|)|x|^{-k}$ is not a positive definite distribution on $\mathbb{R}^n$, where $|\cdot|$ is
the Euclidean norm on $\mathbb{R}^n$. Then there exists an even function $g\in C^\infty (S^{n-1})$
such that

\begin{equation}\label{Zhang1}
\int_{S^{n-1}}f(x) g(x) dx >0
\end{equation} and
\begin{equation}\label{Zhang2}
\int_{S^{n-1}\cap H} g(x) dx\le 0,
\end{equation}
for any $(n-k)$-dimensional plane $H$ through the origin.
\end{Lem}
\proof

Since $f$ is infinitely differentiable, by \cite[Section 3.2]{Kbook}, $(f(x/|x|)|x|^{-k})^\wedge$
is a continuous function on $\mathbb{R}^n\setminus \{0\}$. By our assumption there exists $\xi\in
S^{n-1}$ such that $(f(x/|x|)|x|^{-k})^\wedge(\xi)<0$. By continuity of $(f(x/|x|)|x|^{-k})^\wedge$
there is a neighborhood of $\xi$ where this function is negative. Let
$$\Omega=\{\theta \in S^{n-1}: (f(x/|x|)|x|^{-k})^\wedge(\theta)<0\}.$$
Choose a non-positive infinitely-smooth even function $v$ supported in $\Omega$. Extend $v$ to a
homogeneous function $|x|^{-n+k}v(x/|x|)$ of degree $-n+k$ on $\mathbb{R}^n$. By \cite[Section
3.2]{Kbook}, the Fourier transform of $ |x|^{-n+k}v(x/|x|)$ is equal to $|x|^{-k}g(x/|x|)$ for some
infinitely smooth function $g$ on $S^{n-1}$.

By Parseval's formula on the sphere (Lemma \ref{Lem:Parseval}) we have
\begin{eqnarray*}
&&\hspace{-1cm}\int_{S^{n-1}} f(x) g(x) dx= \int_{S^{n-1}}\left( f(x/|x|) |x|^{-k}\right)
\left(g(x/|x|) |x|^{-n+k}\right)
dx\\
&=&\frac{1}{(2\pi)^n} \int_{S^{n-1}}\left( f(x/|x|) |x|^{-k} \right)^\wedge (\theta) \left(g(x/|x|)
|x|^{-n+k}\right)^\wedge (\theta) d\theta\\
&=& \frac{1}{(2\pi)^n} \int_{S^{n-1}}\left( f(x/|x|) |x|^{-k} \right)^\wedge (\theta) v (\theta)
d\theta >0,
\end{eqnarray*}
since $v$ is non-positive and supported in the set where $\left( f(x/|x|) |x|^{-k} \right)^\wedge$
is negative.

Secondly, by Lemma \ref{Lem:FTonH} we have

\begin{eqnarray*}
&&(2\pi)^k \int_{S^{n-1}\cap H} g(x) dx=(2\pi)^k \int_{S^{n-1}\cap H} g(x/|x|) |x|^{-n+k} dx\\
&&=\int_{S^{n-1}\cap H^\perp} \left(g(x/|x|) |x|^{-n+k}\right)^\wedge (\theta) d\theta=
\int_{S^{n-1}\cap H^\perp} v(\theta) d\theta \le 0,
\end{eqnarray*}
since $v$ is non-positive.

\qed

\section{Main results}

\begin{Prop}\label{NotIntBodies} Let $1\le k \le n-2$.
There exists an infinitely smooth origin symmetric strictly e-convex body $L$ in the unit ball
$B^n\subset\mathbb{R}^n$, so that
\begin{equation}
\frac{\|x\|^{-k}_L}{(1-(\frac{|x|}{\|x\|_L})^2)^k}
\end{equation}
is not a positive definite distribution on $\mathbb{R}^n$.
\end{Prop}
\proof First, we consider the cases $k=n-2$ and $n-3$. We will  use a construction similar to
\cite[Proposition 3.9]{Y}. Let $L$ be a circular cylinder of radius $\sqrt{2}/2$ with  $x_n$ being
its  axis of revolution. To the top and bottom of  the cylinder attach spherical caps, that are
totally geodesic in the spherical metric. Clearly the body $L$ constructed this way is e-convex and
therefore h-convex.
% as a body in $\mathbb{H}^n$ (in the sense that it has nonzero principal curvatures of the same sign).
Using the formula
\begin{equation}\label{eqn:M}
\|x\|^{-1}_M=\displaystyle\frac{\|x\|^{-1}_L}{1-(\frac{|x|}{\|x\|_L})^2}
\end{equation} we define a body $M$. (Note, that $M$ is well-defined, since $L$ lies entirely in the unit ball $B^n$ and the
denominator in the latter formula is never equal to zero).

\begin{center}
\includegraphics{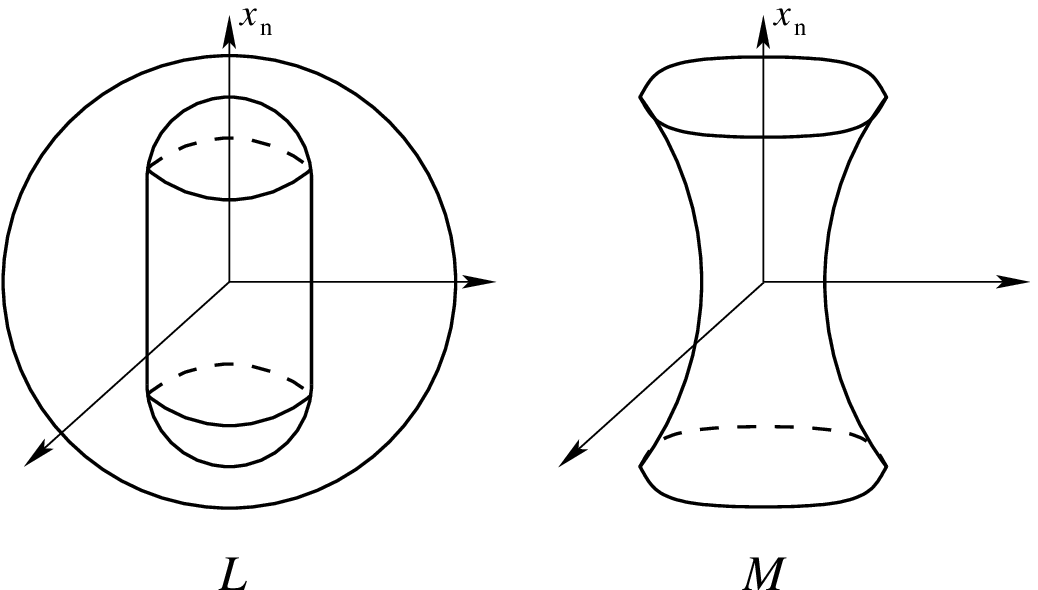}
\nopagebreak
 %Figure 2
\end{center}
Clearly the body $M$ is the image of $L$ under the map:
\begin{eqnarray}\label{eqn:map}
\displaystyle (r,\theta)\mapsto\left(\frac{r}{1-r^2},\theta\right)
\end{eqnarray}
It can be checked directly that the cylinder is mapped into the surface of revolution obtained by
rotating the hyperbola $x_1=\displaystyle \frac{1}{2}\left(\sqrt{2}+\sqrt{2+4x_n^2}\right)$ about
the $x_n$-axis, and the top and bottom spherical caps are mapped into flat disks. The latter
follows from the fact that (\ref{eqn:map}) maps s-geodesics into e-geodesics. Indeed, without loss
of generality we may consider a s-geodesic given by the equation: $r^2+a\ r \cos\phi - 1=0$ in some
2-dimensional plane. The image of this s-geodesic under the map (\ref{eqn:map}) is an e-geodesic
$\displaystyle r=\frac{1}{a\cos\phi}$.

The body $L$ constructed above is not smooth. But we can approximate it by infinitely smooth
e-convex bodies that differ from $L$ only in a small neighborhood of the edges.
%Such approximation affects $A_{M,\xi}(t)$ only in a neighborhood of the edges.
Since the body $M$ is obtained from $L$
by (\ref{eqn:M}), and the denominator in (\ref{eqn:M}) is never equal to zero, the body $M$ is also
infinitely smooth.

%(Now that the bodies $L$ and $M$ are smooth, Figure 2 might be confusing, but we wanted to make it
%as simple as possible, just to emphasize the idea).

Now that we have defined the body $M$, we can explicitly compute its parallel section function
$A_{M,\xi}$ in the direction of the $x_n$-axis.

\begin{eqnarray*}
A_{M,\xi}(t)= C_n \displaystyle  \left({\sqrt{2}+\sqrt{2+4t^2}}\right)^{n-1}.
\end{eqnarray*}

Let the height of the cylindrical part of $L$ be equal to $\sqrt{2}-2\lambda$ and the height of its
image under (\ref{eqn:map}) equal to $2N$ (see the picture below). Since the radius of the cylinder
equals $\sqrt{2}/2$, when $\lambda$ tends to zero the top and bottom parts of the body $L$ get
closer to the sphere $x_1^2+\cdots+x_n^2=1$. Recalling the definition of the radial function of
$M$:
$$\rho_M(x)=\displaystyle\frac{\rho_L(x)}{1-\rho_L(x)^2},\quad \forall x\in S^{n-1},$$
one can see that the height $2N$ of the body $M$ approaches infinity as $\lambda \to 0$.

\begin{center}
\includegraphics{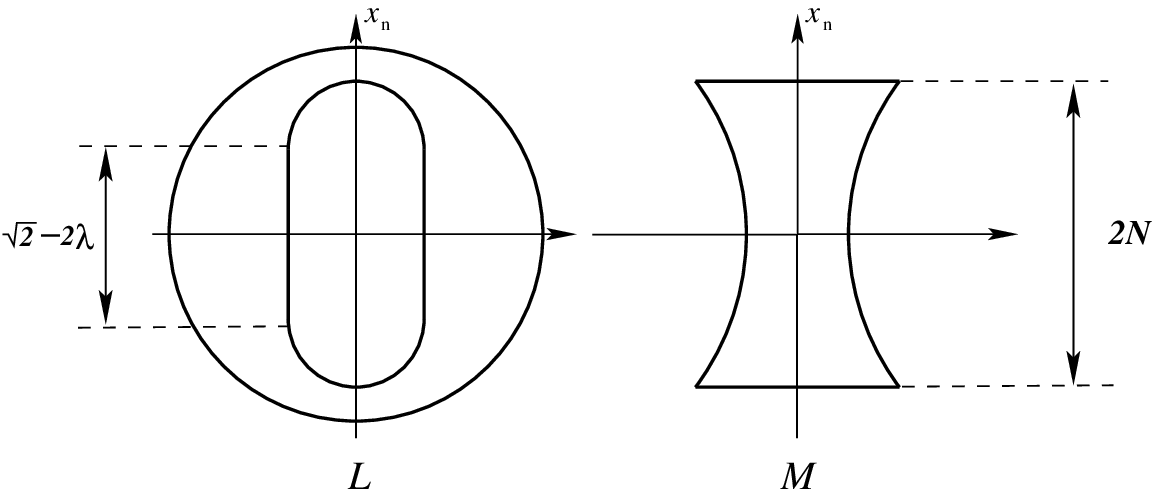}
\nopagebreak
 %Figure 2
\end{center}

Since $M$ is an infinitely smooth body,  $(\|x\|_M^{-n+k+1})^\wedge$ is a function. Applying
Theorem \ref{Thm:GKS}  with $k=1$ we get

\begin{eqnarray*}
(\|x\|_M^{-n+2})^\wedge(\xi)&=&-2(n-2) \int_0^\infty \frac{A_{M,\xi}(t)-A_{M,\xi}(0)}{t^2}dt\\
                            &=&- 2(n-2) C_n\int_0^N \frac{\left({\sqrt{2}+\sqrt{2+4t^2}}\right)^{n-1}-(2\sqrt{2})^{n-1}}{t^2}dt+\\
                            & &+ 2(n-2)C_n
                                \int_N^\infty\frac{(2\sqrt{2})^{n-1}}{t^2}dt.
\end{eqnarray*}

%Consider $\left({\sqrt{2}+\sqrt{2+4t^2}}\right)^{n-1}$,
To estimate the first integral we use the binomial theorem,
\begin{eqnarray*}
\left({\sqrt{2}+\sqrt{2+4t^2}}\right)^{n-1} =(\sqrt{2})^{n-1}+(n-1)(\sqrt{2})^{n-2}\sqrt{2+4t^2}+\\
+\frac{(n-1)(n-2)}{2} (\sqrt{2})^{n-3}({2+4t^2})+\cdots\\
\ge (2\sqrt{2})^{n-1} + 2 {(n-1)(n-2)} (\sqrt{2})^{n-3}t^2,
\end{eqnarray*}
where the last inequality was obtained by putting $t=0$ in all the terms of the binomial expansion,
except for the third term. Therefore, for some positive constants $C'_n$ and $C''_n$ we have
$$
(\|x\|_M^{-n+2})^\wedge(\xi)\le-C'_n\int_0^N dt+C''_n\int_N^\infty\frac{1}{t^2}dt=
                            -C'_n N+C''_n\frac{1}{N}<0
$$
for $N$ large enough.

Therefore the body $M$, corresponding to this $N$, is not a $(n-2)$-intersection body in the
Euclidean sense, which implies that

\begin{equation}
\frac{\|x\|^{-n+2}_L}{(1-(\frac{|x|}{\|x\|_L})^2)^{n-2}}= \|x\|^{-n+2}_M
\end{equation}
is not a positive definite distribution.

Similarly we can show that $M$ is not a $(n-3)$-intersection body. Indeed, if $k=2$ Theorem
\ref{Thm:GKS} implies
$$(\|x\|_M^{-n+3})^\wedge(\xi)=-\pi (n-3) A_{M,\xi}^{''}(0)< 0,$$ since the second derivative of
the function $A_{M,\xi}$ equals:
$$A_{M,\xi}^{''}(0)=C_n(n-1)(2\sqrt{2})^{n-1}>0 .$$

Next we  handle the case when $1\le k < n-3$. For this we use a different construction.  Let $M$ be
an infinitely smooth origin symmetric e-convex body in $\mathbb{R}^n$, for which $\|x\|_M^{-k}$ is
not positive definite. %(For example, the unit balls of the spaces $\ell^n_q$, $2<q\le \infty$, see
%\cite{K1}. In fact not all of them are infinitely smooth, but one can approximate them by
%infinitely smooth bodies, that are not $k$-intersection bodies for $1\le k < n-3$).
(For example, the unit ball of the space $\ell^n_4$, see \cite{K1}). Dilate this body $M$, if
needed, to make sure that it lies in the unit Euclidean ball. Let $\rho_M(x)$ be the radial
function of this body. Define a body $L$ as follows:
$$\rho_L(x)=\frac{-1+\sqrt{1+4(\rho_M(x))^2}}{2\rho_M(x)}, \quad \mbox{ for } x\in S^{n-1}.$$
One can check that
$$\rho_M(x)=\displaystyle\frac{\rho_L(x)}{1-\rho_L(x)^2},\quad \mbox{ for } x\in S^{n-1}.$$
Clearly, $M$ is the image of $L$ under the transformation (\ref{eqn:map}). Since (\ref{eqn:map})
maps s-geodesics into e-geodesics, $L$ is a s-convex body, and therefore e-convex.

%We claim that $L$ is s-convex (and therefore e-convex). Suppose it is not. Then there exist two
%points in $ L$, such that the segment of the s-geodesic connecting them does not lie in $L$.
%Without loss of generality we may assume that this s-geodesic is given by the equation: $r^2+a\ r
%\cos\phi - 1=0$ in some 2-dimensional plane. The image of this s-geodesic under the map
%(\ref{eqn:map}) is an e-geodesic $\displaystyle r=\frac{1}{a\cos\phi}$. So, the image of the
%s-geodesic under the map (\ref{eqn:map}) is a straight line segment, which does not lie entirely in
%$M$. This is a contradiction, since  $M$ is e-convex.

Thus we have proved that for $1\le k <n-3$,

$$\frac{\|x\|^{-k}_L}{\left(1-(\frac{|x|}{\|x\|_L})^2\right)^k}=\|x\|_M^{-k}$$
is not positive definite.

To finish the proof, note that in our construction $L$ is not necessarily strictly e-convex. But
one can replace $L$ with $L_\epsilon$, defined by
$$\|\theta\|_{L_\epsilon}^{-1}= \|\theta\|_{L}^{-1}+\epsilon |\theta|^{-1}.$$
One can choose $\epsilon>0$ small enough, so that $L_\epsilon$ is strictly e-convex, and so that
$\frac{\|x\|^{-k}_{L_\epsilon}}{\left(1-(\frac{|x|}{\|x\|_{L_\epsilon}})^2\right)^k}$ is still not
positive definite (see, for example, the approximation argument in \cite[Lemma 4.10]{Kbook}).
 \qed

\begin{Thm} Let $1\le k <n-1$.
There are origin-symmetric  convex bodies $K$ and $L$ in $\mathbb{H}^n$, $n\ge 3$, such that
$$\mathrm{vol}_{n-k}(K\cap H)\le \mathrm{vol}_{n-k}(L\cap H)$$ for every $(n-k)$-dimensional totally geodesic
plane through the origin, but
$$\mathrm{vol}_{n}(K)> \mathrm{vol}_{n}(L).$$
\end{Thm}

\proof Let $L$ be an infinitely smooth origin symmetric e-convex body from Proposition
\ref{NotIntBodies}, for which
$\displaystyle\frac{||x||^{-k}_L}{\left(1-(\frac{|x|}{||x||_L})^2\right)^k}$ is not positive
definite.

By Lemma \ref{Lem:Zhang} there exists an even function $g\in C^\infty (S^{n-1})$ such that

\begin{equation}\label{Zhang1}
\int_{S^{n-1}}\frac{||x||^{-k}_L}{\left(1-(\frac{|x|}{||x||_L})^2\right)^k} g(x) dx >0
\end{equation} and
\begin{equation}\label{Zhang2}
\int_{S^{n-1}\cap H} g(x) dx\le 0, \quad \mbox{ for all}\,\, H.
\end{equation}

Now apply a standard argument to construct another body $K$ which along with the body $L$ provides
a counterexample to the hyperbolic LDBP problem (cf.  \cite{K3}, Theorem 2 or \cite{Zv}, Theorem
2).
%By continuity of
%$\left(\frac{||x||^{-k}_L}{\left(1-(\frac{|x|}{||x||_L})^2\right)^k}\right)^\wedge$ on the unit
%sphere there is a neighborhood  where this function is negative. Let
%$$\Omega=\{\theta \in S^{n-1}: \left(\frac{||x||^{-k}_L}{\left(1-(\frac{|x|}{||x||_L})^2\right)^k}\right)^\wedge(\theta)<0\}.$$
%Choose a non-positive infinitely-smooth even function $v$ supported on $\Omega$. Extend $v$ to a
%homogeneous function $r^{-n+k}v(\theta)$ of degree $-n+k$ on $\mathbb{R}^n$. By  \cite[Lemma
%5]{K3}, the Fourier transform of $ r^{-n+k}v(\theta)$ is equal to $r^{-k}g(\theta)$ for some
%infinitely smooth function $g$ on $S^{n-1}$.
Define a new body $K$ as follows:
\begin{eqnarray}\label{DefK}
\int_{0}^{\|\theta\|^{-1}_K} \frac{r^{n-k-1}}{(1-r^2)^{n-k}}dr=\int_{0}^{\|\theta\|^{-1}_L}
\frac{r^{n-k-1}}{(1-r^2)^{n-k}}dr+\epsilon g(\theta)
\end{eqnarray}
for $\theta\in S^{n-1}$ and some  $\epsilon>0$  small enough (to guarantee that $K$ is still convex
in hyperbolic sense). Indeed, define a function $\alpha_\epsilon (\theta)$ such that
\begin{eqnarray*}
\int_{0}^{\|\theta\|^{-1}_L} \frac{r^{n-k-1}}{(1-r^2)^{n-k}}dr+\epsilon v(\theta)=
\int_{0}^{\|\theta\|^{-1}_L+\alpha_\epsilon (\theta)} \frac{r^{n-k-1}}{(1-r^2)^{n-k}}dr,
\end{eqnarray*}
then
\begin{eqnarray*}
\|\theta\|^{-1}_K= \|\theta\|^{-1}_L+\alpha_\epsilon (\theta).
\end{eqnarray*}
The function $\alpha_\epsilon (\theta)$ and its first and second derivatives converge uniformly to
zero as $\epsilon\to 0$ (cf. \cite[Proposition 2]{Zv}), therefore since $L$ is strictly e-convex,
there exists $\epsilon$ small enough, so that $K$ is also strictly e-convex, and hence h-convex.

Let $H$ be an $(n-k)$-plane through the origin. Integrating (\ref{DefK}) over  $S^{n-1}\cap H$ and
using inequality (\ref{Zhang2}), we get
\begin{eqnarray*}
\int_{S^{n-1}\cap H}\int_{0}^{\|\theta\|^{-1}_K} \frac{r^{n-k-1}}{(1-r^2)^{n-k}}dr d\theta\le
\int_{S^{n-1}\cap H}\int_{0}^{\|\theta\|^{-1}_L} \frac{r^{n-k-1}}{(1-r^2)^{n-k}}dr d\theta,
\end{eqnarray*}
which, by formula (\ref{eqn:section-polarvolume}), is equivalent to
$$\mathrm{vol}_{n-k}(K\cap H)\le \mathrm{vol}_{n-k}(L\cap H).$$
On the other hand, multiplying both sides of (\ref{DefK}) by
$\left(\frac{\|x\|^{-1}_L}{1-\|x\|^{-2}_L}\right)^k$ and integrating over the sphere $S^{n-1}$ we
get
\begin{eqnarray*}
&&\int_{S^{n-1}}\left(\frac{\|x\|^{-1}_L}{1-\|x\|^{-2}_L}\right)^k\int_{0}^{\|x\|^{-1}_K}
\frac{r^{n-k-1}}{(1- r^2)^{n-k}}drdx= \hspace{4cm}\\
&&=\int_{S^{n-1}}\left(\frac{\|x\|^{-1}_L}{1- \|x\|^{-2}_L}\right)^k\int_{0}^{\|x\|^{-1}_L}
\frac{r^{n-k-1}}{(1- r^2)^{n-k}}drdx +\\
&&\hspace{6cm}  + \epsilon \int_{S^{n-1}}\left(\frac{\|x\|^{-1}_L}{1-\|x\|^{-2}_L}\right)^k g(x)
dx.
\end{eqnarray*}
From   (\ref{Zhang1}) it follows that
\begin{eqnarray*}
&&\int_{S^{n-1}}\left(\frac{\|x\|^{-1}_L}{1-\|x\|^{-2}_L}\right)^k\int_{0}^{\|x\|^{-1}_K}
\frac{r^{n-k-1}}{(1- r^2)^{n-k}}drdx > \hspace{4cm}\\
&&\hspace{3cm}> \int_{S^{n-1}}\left(\frac{\|x\|^{-1}_L}{1-
\|x\|^{-2}_L}\right)^k\int_{0}^{\|x\|^{-1}_L} \frac{r^{n-k-1}}{(1- r^2)^{n-k}}drdx
\end{eqnarray*}

Therefore,
\begin{eqnarray}\label{eqn:beforeZv}
0<\int_{S^{n-1}}\left(\frac{\|x\|^{-1}_L}{1-
\|x\|^{-2}_L}\right)^k\int_{{\|x\|^{-1}_L}}^{\|x\|^{-1}_K} \frac{r^{n-k-1}}{(1- r^2)^{n-k}}drdx
\end{eqnarray}

Next we need the following elementary inequality (cf. Zvavitch, \cite{Zv}). For any $a,b\in (0,1)$
\begin{eqnarray*}
\frac{a^k}{(1- a^2)^k}\int_a^b \frac{r^{n-k-1}}{(1- r^2)^{n-k}}dr\le \int_a^b \frac{r^{n-1}}{(1-
r^2)^{n}}dr.
\end{eqnarray*}
Indeed, since the function $\displaystyle \frac{r^k}{(1-r^2)^k}$ is increasing on the interval
$(0,1)$ we have the following
\begin{eqnarray*}
\frac{a^k}{(1- a^2)^k}\int_a^b \frac{r^{n-k-1}}{(1- r^2)^{n-k}}dr &=&\int_a^b\frac{r^{n-1}}{(1-
r^2)^{n}}\frac{a^k}{(1- a^2)^k}\left(\frac{r^k}{(1-
r^2)^k}\right)^{-1} dr \\
&\le& \int_a^b \frac{r^{n-1}}{(1- r^2)^{n}}dr.
\end{eqnarray*}
Note that in the latter inequality it does not matter whether $a\le b$ or $a\ge b$.

Applying the elementary inequality to (\ref{eqn:beforeZv}) with $a= \|x\|^{-1}_L$ and
$b=\|x\|^{-1}_K$, we get
\begin{eqnarray*}
0&<&\int_{S^{n-1}}\left(\frac{\|x\|^{-1}_L}{1-
\|x\|^{-2}_L}\right)^k\int_{{\|x\|^{-1}_L}}^{\|x\|^{-1}_K} \frac{r^{n-k-1}}{(1-
r^2)^{n-k}}drdx\\
&\le& \int_{S^{n-1}}\int_{{\|x\|^{-1}_L}}^{\|x\|^{-1}_K} \frac{r^{n-1}}{(1- r^2)^{n}}drdx.
\end{eqnarray*}
Hence
\begin{eqnarray*}
\int_{S^{n-1}}\int_0^{{\|x\|^{-1}_L}} \frac{r^{n-1}}{(1- r^2)^{n}}drdx <
\int_{S^{n-1}}\int_{0}^{\|x\|^{-1}_K} \frac{r^{n-1}}{(1- r^2)^{n}}drdx,
\end{eqnarray*}
that is $$\mathrm{vol}_{n}(L)< \mathrm{vol}_{n}(K).$$

 \qed

{\bf Acknowledgments}. The author is thankful to A.Koldobsky for reading this manuscript and making
many valuable suggestions.


\begin{thebibliography}{WWW}

\bibitem[A]{A}
Yu.A.Aminov, {\it The geometry of submanifolds}, Gordon and Breach Science Publishers, Amsterdam,
2001.


\bibitem[BZ]{BZ}
J.Bourgain, Gaoyong Zhang, {\it On a generalization of the Busemann-Petty problem}, Convex
geometric analysis (Berkeley, CA, 1996), 65-76, Math. Sci. Res. Inst. Publ., {\bf 34}, Cambriddge
Univ.Press, Cambridge, 1999.

\bibitem[DFN]{DFN}
B.A.Dubrovin, A.T.Fomenko, S.P.Novikov, {\it Modern geometry - methods and applications. Part I.
The geometry of surfaces, transformation groups, and fields. Second edition}, Springer-Verlag, New
York, 1992.


%\bibitem[GHS]{GHS}
%F.Gao, D.Hug, R.Schneider, {\it Intrinsic volumes and polar sets in spherical space}, Math. Notae,
%A\~{n}o  XLI (2001/02), 159--176 (2003).

%\bibitem[G]{G} R.J.Gardner, {The Brunn-Minkowski inequality}, Bulletin (New Series) of the American
%Mathematical Society, Vol. 39, Number 3, 355-405.

\bibitem[GKS]{GKS}
R.J.Gardner, A.Koldobsky, T.Schlumprecht, {\it An analytic solution to the Busemann-Petty problem
on sections of convex bodies}, Annals of Math. {\bf 149} (1999), 691--703.

\bibitem[GV]{GV}
I.M.Gelfand, N.Ya.Vilenkin, {\it Generalized functions, vol.4. Applications of harmonic analysis},
Academic Press, New York, 1964.


%\bibitem[K]{K}
%A.Koldobsky,  {\it The Busemann-Petty problem via spherical harmonics}, Advances in Math. {\bf 177} (2003), 105--114.

%\bibitem[K1]{K1}
%A.Koldobsky,  {\it Inverse formula for the Blashke-Levy representation}, Houston J. Math. {\bf 23}
%(1997), 95--108.

\bibitem[K1]{K1}
A.Koldobsky, {\it Intersection bodies in $\mathbb{R}^4$}, Advances in Math, {\bf 136} (1998),
1--14.


\bibitem[K2]{K3} A.Koldobsky,  {\it A generalization of the Busemann-Petty problem on sections of
convex bodies}, Israel J. Math. {\bf 110} (1999), 75--91.

\bibitem[K3]{K4}
A.Koldobsky, {\it A functional analytic approach to intersection bodies}, Geom. Funct. Anal. {\bf
10}(2000),1507-1526.

\bibitem[K4]{Kbook}
A.Koldobsky, {\it Fourier analysis in convex geometry}, to appear.


%\bibitem[K5]{K5}
%A.Koldobsky, {\it On the derivatives of x-ray functions}, Arch. Math. {\bf 79} (2002), 216-222.

%\bibitem[K6]{K6}
%A.Koldobsky,  {\it The Busemann-Petty problem via spherical harmonics}, Advances in Math. {\bf 177}
%(2003), 105--114.


%\bibitem[KYY]{KYY}
%A.Koldobsky, V.Yaskin, M.Yaskina, {\it Modified Busemann-Petty problem on sections  of convex
%bodies}, preprint.

%\bibitem[L]{L}
%E.Lutwak, {\it Intersection bodies and dual mixed volumes}, Advances in Math., {\bf 71} (1988),
%232--261.

\bibitem[MP]{MP}
D.Mej\'{i}a, Ch.Pommerenke, {\it On spherically convex univalent functions}, Michigan Math. J.,
{\bf 47} (2000), 163--172.


\bibitem[P]{P}
A.V.Pogorelov, {\it Extrinsic geometry of convex surfaces}, Translations of Mathematical
Monographs, vol.35, American Mathematical Society, Providence, RI, 1973.


\bibitem[R]{R}
J.G.Ratcliffe, {\it Foundations of hyperbolic manifolds}, Springer-Verlag, New York, 1994.

%\bibitem[RZ]{RZ}
%B.Rubin, Gaoyong Zhang, {\it Generalizations of the Busemann-Petty problem for sections of convex
%bodies}, J. Funct. Anal., {\bf 213} (2004), 473--501.

\bibitem[Y]{Y}
V.Yaskin, {\it The Busemann-Petty problem in hyperbolic and spherical spaces}, preprint.


\bibitem[Zh1]{Zh1}
Gaoyong Zhang, {\it Sections of convex bodies }, Amer. J. Math. {\bf  118} (1996), 319--340.

\bibitem[Zh2]{Zh2}
Gaoyong Zhang, {\it A positive answer to the Busemann-Petty problem in four dimensions}, Annals of
Math. {\bf  149} (1999), 535-543.

\bibitem[Zv]{Zv}
A.Zvavitch, {\it The Busemann-Petty problem for arbitrary measures}, Math. Ann., to appear.

\end{thebibliography}
\end{document}